\documentclass[a4paper,oneside,12pt]{article}
\usepackage{amssymb}

\newtheorem{theorem2}{Th\'eor\`eme}[section]

\newtheorem{lemma2}{Lemme}[section]

\newtheorem{proposition2}{Proposition}[section]

\newtheorem{remark2}[theorem2]{Remarque \rm}

\newcommand{\ds}{\displaystyle}
\newcommand{\ofr}{{\mathfrak o}}
\newcommand{\pfr}{{\mathfrak p}}

\newcommand{\lra}{\longrightarrow}
\newcommand{\noi}{\noindent}

\newcommand{\CC}{\mathbb C}

\newcommand{\VV}{\mathcal V}

\newcommand{\Afr}{\mathfrak A}
\newcommand{\Bfr}{\mathfrak B}
\newcommand{\Cfr}{\mathfrak C}

\newcommand{\Pfr}{\mathfrak P}

\title{Formules de caract\`ere pour la s\'erie discr\`ete de ${\rm
    GL}(N)$ \\
Guide de l'utilisateur pour  typistes}
\author{Paul Broussous\\
Universit\'e de Poitiers}

\begin{document}
\maketitle

\noi {\Large Introduction --  } Ces notes s'adressent \`a des
sp\'ecialistes de la th\'eorie des types pour les groupes r\'eductifs
$p$-adiques et je m'excuse par avance aupr\`es du lecteur qui
trouvera tr\`es peu de rappels dans ce qui suit, autant en ce qui
concerne la th\'eorie que les notations. 
\medskip

 L'explicitation de la correspondance
de Jacquet-Langlands a \'et\'e faite par Silberger et Zink en niveau
$0$, et par Bushnell et Henniart,  en niveau quelconque, pour un grand
nombre de repr\'esentations supercuspidales de ${\rm GL}(N,F)$. 
L'extension des travaux de Bushnell et Henniart au cas des
repr\'esentations non supercuspidales  de la s\'erie discr\`ete
demande d'\'etablir des formules pour la valeur du caract\`ere
d'Harish-Chandra en les \'el\'ements elliptiques
r\'eguliers. Encore r\'ecemment peu de formules \'etaient
connues. 
\smallskip

J'ai pu obtenir dans [Br] quelques premi\`ere formules
 particuli\`eres. L'id\'ee de base, due \`a Henniart,  \'etait de
 transf\'erer le pseudo-coefficient de Kottwitz par les isomorphismes
 d'alg\`ebres de Hecke de [BK].
\smallskip

 Dans le travail [BS] que je pr\'esente ici (en commun avec P. Schneider),
 je propose des formules tout-\`a-fait g\'en\'erales, retrouvant comme
 cas particuliers les formules obtenues dans [Br]. Cette fois-ci,
 l'id\'ee n'est plus de transf\'erer le pseudo-coefficient de
 Kottwitz, mais d'en construire un directement pour tout membre de la
 s\'erie discr\`ete, par des m\'ethodes de nature homologique (tout
 comme chez Kottwitz).
\smallskip

 Ces pseudo-coefficients s'obtiennent en construisant des syst\`emes
 de coefficients \'equivariants sur l'immeuble. Une premi\`ere
 construction  [SS] avait \'et\'e faite par
 Schneider et Stuhler dans le cas d'un groupe r\'eductif
 quelconque. Cependant, sauf en niveau $0$, ces syst\`emes de
 coefficients ne donnaient pas lieu \`a des formules de caract\`ere
 exploitables. 
\smallskip

 Avec Schneider, nous modifions la construction de [SS] en nous
 servant des types simples de Bushnell et Kutzko pour construire des
 syst\`emes de coefficients. C'est Schneider le
 premier qui a devin\'e que, cach\'ee dans la monographie [BK],
 il y a de l'homologie sur l'immeuble de Bruhat-Tits. 
\smallskip

 L'objet de ces notes est de fournir \`a un utilisateur potentiel de
 nos  formules de caract\`ere un \'enonc\'e succint, mais complet et
 utilisable, des r\'esultats de
 [BS]. Elles s'adressent \`a des personnes suffisamment \`a l'aise
 avec les notations et techniques de [BK].

\tableofcontents

\section{Notations : groupes et immeubles}

 Pour les assertions non prouv\'ees et/ou non r\'ef\'erenc\'ees, on
 renvoie le lecteur \`a la monographie de Bushnell et Kutzko
 [BK], ou bien au {\S}I de [BS].
\medskip

 Si $K$ est un corps localement compact et non archim\'edien, on note
 $\ofr_F$ son anneau d'entiers, $\pfr_K$ l'id\'eal maximal de $\ofr_K$
 et $k_K = \ofr_K /\pfr_K$ le corps r\'esiduel (fini). On fixe une
 fois pour toute un tel corps $F$. 
\smallskip

 Soit $V$ un $F$-espace vectoriel de dimension finie $N$. On pose
 $A={\rm End}_F (V)\simeq {\rm M}(N,F)$. Si $E/F$ est une extension de
 corps finie de $F$ plong\'ee dans $A$, le commutant de $E$ dans $A$
 est $B = {\rm End}_E (V)\simeq {\rm M}(N/[E:F], E)$.  On note $G$ le
 groupe ${\rm Aut}_F (V)$ et $G_E$ le groupe ${\rm Aut}_E (V)$,
 centralisateur de $E^\times$ dans $G$. 

 On note ${\rm Her}(A)$ (resp ${\rm Her}(B)$) l'ensemble des
 $\ofr_F$-ordres h\'er\'editaires dans $A$ (resp. des $\ofr_E$-ordres
 h\'er\'editaires dans $B$). Ce sont des ensembles partiellement
 ordonn\'es (EPO) munis des actions par conjugaison de $G$ et $G_E$
 respectivement. 

 On a une injection  $G_E$-\'equivariante $j_{E/F}$~: ${\rm
   Her}(A)\lra {\rm Her}(B)$. Elle associe \`a un ordre
 h\'er\'editaire $\Bfr$ de $B$ attach\'e \`a une $\ofr_E$-cha\^ine de
 r\'eseaux $({\mathcal L})$ de $V$, l'ordre h\'er\'editaire $\Afr =
 \Afr (\Bfr )$ associ\'e \`a $({\mathcal L})$ vue comme
 $\ofr_F$-cha\^ine de r\'eseaux.   

L'immeuble semisimple $X_F$ de $G$ est naturellement  la r\'ealisation
g\'eom\'etrique  d'un  complexe simplicial de
dimension $N-1$, et de  fa\c con abusive, on notera encore $X_F$  ce
complexe simplicial. C'est un espace topologique localement compact sur
lequel $G_F$ agit par automorphismes simpliciaux. Il existe une
bijection d\'ecroissante et $G_F$-\'equivariante entre l'EPO ${\rm
  Her}(A)$ et l'EPO des simplexes de $X_F$, que l'on notera $\Afr
\mapsto \sigma_\Afr$. Elle est caract\'eris\'ee comme suit : si $\Afr$
est un ordre h\'er\'editaire de $A$, $\sigma_\Afr$ est l'unique
simplexe de l'immeuble dont le fixateur compact est le sous-groupe
parahorique  $U(\Afr
)=\Afr^\times$.  On a des notations et faits similaires pour le groupe
$G_E$. 
\medskip

Sur les r\'ealisations g\'eom\'etriques $X_F$ et $X_E$, on a
  des structures affines : le barycentre de deux points \`a
  coefficients positifs est d\'efini. Le fait suivant sera un
  ingr\'edient crucial de nos constructions.

\begin{theorem2} [BL] i)  Il existe une 
  unique application $j$~: $X_E \lra X_F$ qui est $G_E$-\'equivariante
  et affine. 

ii) De plus les complexes de drapeaux   ${\rm Flag} ({\rm Her}(A))$ et
${\rm Flag}({\rm Her}(B))$ correspondent respectivement aux premi\`eres
  subdivisions barycentriques de $X_F$ et $X_E$. On montre que $j$
  correspond \`a  l'application ${\rm Flag}({\rm Her}(B)) \lra {\rm Flag} ({\rm
    Her}(A))$ induite par $j_{E/F}$. En particulier, si $\Bfr\in {\rm
    Her}(B)$, alors $j$ envoie l'isobarycentre de $\sigma_\Bfr$ sur
  l'isobarycentre de $\sigma_{\Afr (\Bfr )}$.
\end{theorem2}

 Dans la suite on identifiera syst\'ematiquement le $G_E$-ensemble 
 $X_E$ avec son image $j(X_E)$ dans $X_F$.  L'inclusion $X_E \subset
 X_F$ n'est pas simpliciale en g\'en\'eral. Elle l'est si, et
 seulement si, l'extension $E/F$ est non ramifi\'ee. Cependant
 l'inclusion $X_E \subset X_F$ est toujours simpliciale apr\`es
 passage \`a la premi\`ere subdivision barycentrique. 

\medskip

 Le sous-$G_F$-ensemble  $X(E) = \ds \bigcup_{g\in G} g.X_E$ de $X_F$
 est naturellement la r\'ealisation g\'eom\'etrique d'un complexe
 simplicial $X[E]$. Si $E/F$ est non ramifi\'ee, alors $X[E]$ est un
 sous-complexe simplicial de $X_F$. En g\'en\'eral il faut passer
 \`a la premi\`ere subdivision barycentrique pour que l'inclusion
 $X(E)\subset X_F$ soit simpliciale. 
\medskip

 Si $\Afr \in {\rm Her} (A)$, on note $U(\Afr ) = \Afr^\times$ le
 sous-groupe parahorique qui fixe le simplexe $\sigma_\Afr$ de
 $X_F$. Son sous-groupe pro-unipotent est  $U^1 (\Afr ) =
 1+\Pfr_{\Afr}$, o\`u $\Pfr_\Afr$ est le radical de Jacobson de
 $\Afr$.  

 Si $E/F$ est un sous-corps de $A$ et si $\Bfr\in {\rm Her}(B)$, $\Afr
 = \Afr (\Bfr )$, alors on a :
$$
{\mathcal N}(\Afr )\cap G_E = {\mathcal N}(\Bfr ), \ U(\Afr )\cap G_E
= U(\Bfr ), \ U^1 (\Afr )\cap G_E = U^1 (\Bfr ), \ {\mathcal N}(\Afr )
= {\mathcal N}(\Bfr ) U(\Afr )\ .
$$

 De plus, l'action de $G$ sur $X_E$ poss\`ede la propri\'et\'e suivante :
\medskip

\begin{lemma2} ([BS] Lemma I.3.3) Si deux simplexes de $X_E$
  sont conjug\'es sous l'action de $G$,  ils le sont sous l'action de
  $G_E$.
\end{lemma2}

\noi En d'autres termes, $G$ n'induit pas plus d'action sur $X_E$ que $G_E$ le fait d\'ej\`a.

\section{Caract\`eres simples}

  Les r\'ef\'erences pour cette section sont [BK] et
  [BH]. On fixe une fois pour toute une paire simple $[0,\beta ]$
  ([BH](1.5)), c'est-\`a-dire une extension finie $E/F$ munie
  d'un g\'en\'erateur $\beta$ (i.e. $E=F[\beta ]$), satisfaisant les conditions
  suivantes :
\smallskip

[PS1] $\beta\not\in \ofr_E$,

[PS2] $k_0 (\beta , \Afr (E))<0$ (cf. [BK]{\S}1). 
\smallskip

Pour chaque $E$-espace vectoriel $V$ et chaque $\Bfr \in {\rm
  Her}(B)$, o\`u $B={\rm End}_E (V)$, on a une strate simple $[\Afr (\Bfr
  ), n_\Bfr ,0 ,\beta ]$ dans $A={\rm End}_F (V)$, r\'ealisation de
$[0,\beta ]$ dans $A$. 

 Attach\'ees \`a $[\Afr (\Bfr ), n_\Bfr ,0 ,\beta ]$ (donc \`a
 $[0,\beta ]$, $V$, $\Bfr$), on a les donn\'ees suivantes :
\smallskip

 $\bullet$ Deux sous-groupes ouverts compacts de $G={\rm Aut}_F (V)$ : 
$H^1 (\Bfr ) \subset J^1 (\Bfr )\subset U^1 (\Afr (\Bfr ))$, tous deux
normalis\'es par ${\mathcal N}(\Bfr )$. 
\smallskip

 $\bullet$ Un ensemble fini de caract\`eres simples ${\rm \mathcal
  C}(\Bfr ) = {\mathcal C}(\Afr (\Bfr ), 0,\beta )$ de $H^1 (\Bfr )$,
qui ont chacun un $G$-entrelacement donn\'e par $J^1 (\Bfr ) G_E J^1
(\Bfr )$. 
\smallskip

 Rappelons que si $\theta \in {\mathcal C}(\Bfr )$, il existe \`a
 isomorphisme pr\`es une unique rep\'esentation irr\'eductible $\eta =
 \eta (\theta )$ de $J^1 (\Bfr )$ qui contient $\theta$ par
 restriction (la repr\'esentation de Heisenberg de $\theta$). 
\smallskip

 La paire simple et le $E$-espace vectoriel $V$ \'etant fix\'es, on a
 des bijections canoniques :

$$
\tau_{\Bfr_1 ,\Bfr_2}~: \ {\mathcal C}(\Bfr_1 )\lra {\mathcal C}
(\Bfr_2 )\ , \ \Bfr_1 ,\ \Bfr_2 \in  {\rm Her}(B),
$$

\noi appel\'ees {\it applications de transfert}
([BK](3.6)). Gr\^ace \`a ces applications, si l'on fixe une paire
simple $[0,\beta ]$, un $E$-espace vectoriel $V$,  un ordre $\Bfr_0$
dans $B={\rm End}_E (V)$ et un caract\`ere simple $\theta_0\in
{\mathcal C}(\Bfr_0 )$, on obtient une famille de caract\`eres
simples $(H^1 (\Bfr ), \theta (\Bfr ))_{\Bfr \in {\rm Her}(B)}$ en
  posant $\theta (\Bfr ) = \tau_{\Bfr_0 ,\Bfr}(\theta_0 )$. Il lui
 est associ\'e une famille de repr\'esentations de Heisenberg
 $(H^1 (\Bfr ), \theta (\Bfr ))_{\Bfr \in {\rm Her}(B)}$ d\'efinies
 \`a isorphisme pr\`es. Il est facile de v\'erifier que ces deux
 familles sont $G_E$-\'equivariantes en un sens \'evident.

 Pour chaque paire d'ordres h\'er\'editaires $\Bfr_1 \subset \Bfr_2$
 dans ${\rm Her}(B)$, on peut former le groupe $J^1 (\Bfr_1 ,\Bfr_2 )
 = U^1 (\Bfr_1 )J^1 (\Bfr_2 )$. 

\begin{proposition2} ([BK] (5.1.14-16), (5.1.18), (5.1.19))
  Fixons une paire simple $[0,\beta ]$, un $E$-espace vectoriel $V$,
  un ordre $\Bfr_0\in {\rm Her}(B)$ et un caract\`ere simple
  $\theta_0\in {\mathcal C}(\Bfr_0 )$. Alors il existe une unique
  famille de repr\'esentations $\ds\left( J^1 (\Bfr_1 ,\Bfr_2 ),\eta
  (\Bfr_1 ,\Bfr_2 )\right)_{\Bfr_1 \subset \Bfr_2}$ (d\'efinies \`a
  isomorphismes pr\`es) qui \'etend la famille $(J^1 (\Bfr ), \eta
  (\Bfr ))_{\Bfr}$ au sens suivant :
\smallskip

(i) $\eta (\Bfr ,\Bfr ) =\eta (\Bfr )$, $\Bfr\in {\rm Her}(B)$ ;

(ii) $\eta (\Bfr_1 ,\Bfr_2 )_{\vert J^1 (\Bfr_2 )} \simeq \eta (\Bfr_1
)$, $\Bfr_1 \subset \Bfr_2 \in {\rm Her}(B)$ ; 

(iii) les induites suivantes sont irr\'eductibles et \'equivalentes :
$$
{\rm Ind}_{J^1 (\Bfr_1 )}^{U^1 (\Afr_1 )}\, \eta (\Bfr_1 )\simeq {\rm
  Ind}_{J^1 (\Bfr_1 , \Bfr_2 )}^{U^1 (\Afr_1 )}\, \eta (\Bfr_1 ,\Bfr_2
), \    \Bfr_1 \subset \Bfr_2 \in {\rm Her}(B)\ .
$$

De plus on a la relation de compatibilit\'e :
$$
\eta (\Bfr_1 ,\Bfr_2 )_{\vert J^1 (\Bfr_2 ,\Bfr_3)} \simeq \eta
(\Bfr_2 ,\Bfr_3 ), \ \Bfr_1 \subset \Bfr_2 \subset \Bfr_3 \in {\rm
  Her}(B)\ .
$$
\end{proposition2}
\noi Il est facile de v\'erifier que la famille $\ds\left( J^1 (\Bfr_1 ,\Bfr_2 ),\eta
  (\Bfr_1 ,\Bfr_2 )\right)_{\Bfr_1 \subset \Bfr_2}$ de
repr\'esentations de Heisenberg est $G_E$-\'equivariante en un sens
\'evident.

\section{Repr\'esentations de la s\'erie discr\`ete et types simples}

 Sauf si \c ca n'est pas dit express\'ement, on utilisera les
 notations de [BK]. 
\medskip

 On fixe une fois pour toute un type simple $(J,\lambda )$ au sens de
 [BK](5.5.10). On supposera de plus que l'on est en niveau $>0$. 
Concr\`etement cela signifie la chose suivante.
\smallskip

 Il existe une paire simple $[0,\beta ]$, un $E$-espace vectoriel $V$
 de dimension finie, o\`u $E=F[\beta ]$ et un ordre principal $\Bfr_0$
 dans $B ={\rm End}_E (V)$. La repr\'esentation $\lambda$ du groupe $J
 = J(\Bfr_0 ) = J^1 (\Bfr_0 )U (\Bfr_0 )$ est de la forme $\kappa_0
 \otimes \rho$, o\`u $\kappa_0$ est une $\beta$-extension d'un
 caract\`ere simple $\theta_0\in {\mathcal C}(\Bfr_0 )$
 ([BK](5.2.1)), et $\rho$ est l'inflation d'un repr\'esentation
 irr\'eductible cuspidale de $J/J^1 (\Bfr_0 )$ de la forme
 suivante. Rappelons que le quotient  $J/J^1 (\Bfr_0 )$ s'identifie
 \`a ${\rm GL}(n/e, k_E )^{\times e}$, o\`u $n:={\rm dim}_E (V)$, $e$
   est la p\'eriode de l'ordre $\Bfr_0$, et $k_E$ d\'esigne le corps
   r\'esiduel de $E$. Alors la condition sur $\rho$ et que, comme
   repr\'esentation de  ${\rm GL}(n/e, k_E )^{\times e}$, elle est de
     la forme $\rho_0^{\otimes e}$, o\`u $\rho_0$ est une
     repr\'esentation (irr\'eductible, cuspidale) de ${\rm GL}(n/e,k_E
     )$. 
\medskip

 Les donn\'ees d\'ecrivant le type simple $(J,\lambda )$ de $G = {\rm
   Aut}_F (V)$ ne sont pas uniques. C'est pour cette raison que {\it nous
 fixons une fois pour toute} :
\smallskip
  
 --  une paire simple $[0,\beta ]$, 

 -- un $E$-espace vectoriel $V$

 -- un ordre principal $\Bfr_0$ dans $B= {\rm End}_E (V)$,

 -- un caract\`ere simple $\theta_0 \in {\mathcal C}(\Bfr_0 )$,

 -- une $\beta$-extension $\kappa_0$ de $\theta_0$,

 -- une repr\'esentation cuspidale irr\'eductible $\rho_0$ de ${\rm
   GL}(n/e,k_E )$.
\smallskip

 De plus d'apr\`es la Proposition 2.1, ces donn\'ees donnent lieu \`a
 une famille $G_E$-\'equivariante de repr\'esentations de Heisenberg 
 $\ds\left( J^1 (\Bfr_1 ,\Bfr_2 ),\eta  (\Bfr_1 ,\Bfr_2
 )\right)_{\Bfr_1 \subset \Bfr_2}$. 
\smallskip

 Fixons une extension non ramifi\'ee $L/E$ contenue dans $B$, v\'erifiant 
$[L:F]=n/e$, telle que   le groupe multiplicatif $L^\times$ normalise ${\Bfr}_0$.
 On pose ${\rm End}_L \, V \simeq {\rm M}(e,L)$ et 
$G_L = {\rm Aut}_L \, V$. On a une application canonique ${\rm Her}(C)\lra
{\rm Her}(B)$ ainsi qu'une inclusion simpliciale et $G_L$-\'equivariante
$X_L\subseteq X_E$. Notons que l'unique ordre  ${\Cfr }_0 \in {\rm Her}(C)$
v\'erifiant ${\Bfr }_0 ={\Bfr }({\Cfr }_0 )$ est un ordre minimal
(ou {\it ordre d'Iwahori}) ; par contre ${\Bfr}_0$ n'est pas minimal en 
g\'en\'eral. On note $X[L] = \ds \bigcup_{g\in G} gX_L$ que l'on munit de la structure 
simpliciale naturelle $G$-invariant qui prolonge celle de $X_L$. 
\medskip

 On fixe un ordre maximal ${\Cfr }_{\rm max}\supseteq {\Cfr}_0$ et on pose 
${\Bfr}_{\rm max} = {\Bfr}( {\Cfr}_{\rm max})$ ; cet ordre est 
toujours maximal. D'apr\`es [BK](5.2.2-5), à isomorphisme pr\`es, il existe une
unique $\beta$-extension $\kappa_{\rm max}$ de $\eta ({\Bfr}_{\rm max} )$ telle
que
$$
{\rm Ind}_{J({\Bfr}_0 )}^{U({\Bfr}_0 )U^1 ({\Afr}({\Bfr}_0 )}
\kappa_0 \simeq {\rm Ind}_{U({\Bfr}_0)J^1 ({\Bfr}_{\rm
    max})}^{U({\Bfr}_0 )  U^1 ({\Afr}({\Bfr}_0 ))} \kappa_{\rm max}\ .
$$

\noi On peut alors former la repr\'esentation $\lambda_{\rm max} =
 \kappa_{\rm max} \otimes \rho$ de $J_{\rm max} := U(\Bfr_0 )J^1 (\Bfr_{\rm max})$ ; elle 
est irr\'eductible. 

\begin{theorem2} a) La paire $(J_{\rm max},\lambda_{\rm max})$ est un type
de $G$ qui d\'efinit la m\^eme composante de Bernstein
 ${\mathcal R}_{(J,\lambda )}(G)$ que $(J,\lambda )$.  

\noi b) Soit $(\pi ,{\mathcal V})$ un objet de  ${\mathcal R}_{(J,\lambda )}(G)$. On a
${\mathcal V}^{\lambda_{\rm max}} = {\mathcal V}^{\eta (\Bfr_0 ,\Bfr_{\rm max})}$.

\noi c) Soit $(\pi ,{\mathcal V})$ une repr\'esentation irr\'eductible essentiellement 
de carr\'e int\'egrable modulo le centre, objet de  ${\mathcal R}_{(J,\lambda )}(G)$. 
Alors $(\pi ,{\mathcal V})$ contient la paire $(J_{\rm max} ,\lambda_{\rm max})$
(resp. la paire $(J,\lambda )$) avec multiplicit\'e $1$. 
\end{theorem2}

\section{Un syst\`eme de coefficients}

Un syst\`eme de coefficients \'equivariant sur le $G$-complexe simplicial 
$X[L]$ est la donn\'ee de $(({\mathcal V}_{\sigma})_\sigma , (r_\tau^\sigma )_{\tau \subset 
\sigma}, (\varphi_{g,\sigma})_{g,\sigma})$, o\`u :
\medskip

 -- pour chaque simplexe $\sigma$ de $X[L]$, $\VV_\sigma$ est un $\CC$-espace vectoriel,

 -- pour $\tau\subseteq \sigma$, $r_\tau^\sigma\in {\rm Hom}_\CC (\VV_\sigma , 
\VV_{\tau})$. 

 -- pour $\sigma$ simplexe, $g\in G$, $\varphi_{g,\sigma}\in
{\rm Hom}_\CC (\VV_\sigma , \VV_{g\sigma})$. 

 -- pour chaque simplexe $\sigma$, $r_\sigma^\sigma =\varphi_{1,\sigma}={\rm id}_{\VV_\sigma}$,

 -- tous les diagrammes que l'on peut imaginer commutent, 

 -- pour tout $\sigma$ la repr\'esentation  de $G_\sigma$ dans $\VV_\sigma$ induite
par le syst\`eme de coefficients  est lisse.
\medskip

 Soit \`a pr\'esent $(\pi ,\VV )$ une repr\'esentation lisse de $G$. Soit 
$\Cfr \in {\rm Her}(C)$ tel que $\Cfr_{\rm min}\subseteq \Cfr \subseteq \Cfr_{\rm max}$, et
soit $\sigma_\Cfr$ le simplexe de $X_L$ attach\'e \`a $\Cfr$. On pose alors
$$
\VV_{\sigma_\Cfr} = \sum_{g\in U(\Afr )/U(\Bfr )J^1 (\Bfr_{\rm max})} 
\pi (g)\VV^{\eta (\Bfr ,\Bfr_{\rm max})}
$$
\noi o\`u $\Bfr =\Bfr (\Cfr )$ et $\Afr =\Afr (\Bfr )$.

 Si $\sigma$ est un simplexe quelconque de $X[L]$, on peut toujours l'\'ecrire
$\sigma = g.\sigma_{\Cfr}$, pour un $g\in G$ et un $\Cfr\in {\rm Her}(C)$ tel que
$\Cfr_{\rm min}\subseteq \Cfr \subseteq \Cfr_{\rm max}$, et on pose $\VV_\sigma = 
\pi (g) \VV_{\sigma_\Cfr}$.

\begin{theorem2} a) Si $\sigma$ est un simplexe de $X[L]$, $\VV_\sigma$ est bien d\'efini, 
c'est-\`a-dire ne d\'epend d'aucun choix. 

\noi b) Si $\tau\subseteq \sigma$ sont des simplexes de $X[L]$, on a 
$\VV_\sigma\subseteq \VV_\tau$.

\noi c) Si $\sigma$ est un simplexe de $X[L]$ et si $g\in G$, alors $\VV_{g\sigma}
=\pi (g)\VV_\sigma$. 
\end{theorem2}

On peut donc d\'efinir un syst\`eme de coefficients ${\mathcal C}(\pi ) =
 ((\VV_\sigma )_\sigma , (r_\tau^\sigma )_{\tau \subseteq \sigma}, (\varphi_{g,\sigma})_{g,\sigma})$ 
sur $X[L]$, en d\'efinissant $r_\tau^\sigma$ comme \'etant l'inclusion $\VV_\sigma \subseteq \VV_\tau$, et 
$\varphi_{g,\sigma}$ comme \'etant l'application $\VV_\sigma \lra \VV_{g\sigma}$ induite par $\pi (g)$.

\begin{theorem2} Supposons que $(\pi ,\VV )\in {\mathcal R}_{(J,\lambda )}(G)$. Alors le complexe $X[L]$
et le syst\`eme de coefficients ${\mathcal C}(\pi )$ ne d\'ependent que de l'endo-classe $\Theta$
du caract\`ere simple $\theta_0$, et donc d'aucun autre choix fait dans la construction. On  notera
${\mathcal C}_\Theta (\pi )$  ce syst\`eme de coefficients canoniquement attach\'e \`a $\pi$.
\end{theorem2}

 Ce syst\`eme de coefficients peut se calculer presque enti\`erement si $(\pi ,\VV)$ est irr\'eductible 
 et essentiellement de carr\'e int\'egrable, ce que nous supposerons jusqu'\`a la fin de cette section. 
\medskip

Soit $\Cfr\in {\rm Her}(C)$ tel que $\Cfr_{\rm min}\subseteq \Cfr \subseteq 
\Cfr_{\rm max}$. On pose comme d'habitude $\Bfr =\Bfr (\Cfr )$ et 
$\Afr =\Afr (\Bfr )$. Le quotient ${\mathbb G}_\Bfr =U(\Bfr )/U^1 (\Bfr )$
est un produit de groupes  lin\'eaires g\'en\'eraux sur $k_E$, le corps
r\'esiduel de $E$. Il poss\`ede ${\mathbb P}_{\Bfr_0}=U(\Bfr_0)
/U^1 (\Bfr )$ comme sous-groupe parabolique. Ce dernier admet comme
radical unipotent ${\mathbb U}_{\Bfr_0}=U^1 (\Bfr_0)/U^1 (\Bfr )$, et
comme facteur de Levi 
$$
{\mathbb L}_{\Bfr_0} = {\mathbb P}_{\Bfr_0} / {\mathbb U}_{\Bfr_0}
= U(\Bfr_0)/U^1 (\Bfr_0) \simeq {\rm GL}(n/e,k_E )^{\times e}\ .
$$

\noi Soit ${\rm St} (\Bfr ,\rho )$ la repr\'esentation de Steinberg
g\'en\'eralis\'ee de ${\mathbb G}_\Bfr$ de support cuspidal
$({\mathbb L}_{\Bfr_0}, \rho )$ : c'est l'unique sous-repr\'esentation
irr\'eductible g\'en\'erique de l'induite parabolique 
${\rm ind}_{{\mathbb P}_{\Bfr_0}}^{{\mathbb G}_\Bfr}\, \rho$.  On peut alors
former la repr\'esentation $\kappa_{\rm max}\otimes {\rm St}(\Bfr ,\rho )$
de $U(\Bfr )J^1 (\Bfr_{\rm max})$. On montre que cette repr\'esentation est irr\'eductible.

\begin{proposition2}  a) On a $\VV^{\eta (\Bfr ,\Bfr_{\rm max})} \simeq \kappa_{\rm max}\otimes
{\rm St}(\Bfr ,\rho )$ comme $U(\Bfr )J^1 (\Bfr_{\rm max})$-modules. 

\noi b) La repr\'esentation induite :
$$
\lambda (\Afr ):= {\rm ind}_{U(\Bfr )J^1 (\Bfr_{\rm max})}^{U(\Afr )} \, \left(  \kappa_{\rm max}\otimes
{\rm St}(\Bfr ,\rho )\right)
$$
\noi est irr\'eductible.

\noi c) Soit $(\pi ,\VV )\in {\mathcal R}_{(J,\lambda )}(G)$
une repr\'esentation irr\'eductible essentiellement de carr\'e int\'egrable.
Si ${\mathcal C}_\Theta (\pi ) =  
((\VV_\sigma )_\sigma , (r_\tau^\sigma )_{\tau \subseteq \sigma}, (\varphi_{g,\sigma})_{g,\sigma})$ 
est le syst\`eme de coefficients attach\'e \`a $\pi$, alors
$$
\VV_{\sigma_\Cfr} = \VV^{\lambda (\Afr )} \simeq \lambda (\Afr )
$$
\noi o\`u l'isomorphisme est un isomorphisme de $U(\Afr )$-modules.
\end{proposition2}

\begin{remark2} Il est tr\`es important de noter que la donn\'ee du type $(J,\lambda )$ de $\pi$ ne permet pas
de conna\^itre enti\`erement le syst\`eme de coefficients ${\mathcal C}_\Theta (\pi )$, mais seulement les
espaces $\VV_{\sigma_\Cfr }$ comme $U(\Afr )$-modules. Le stabilisateur $G_{\sigma_\Cfr}$ de $\sigma_\Cfr$ dans
$G$  a la forme d'un produit semi-direct $\langle \Pi_\Afr \rangle \ltimes U(\Afr )$, pour un certain 
\'el\'ement $\Pi_\Afr$ de $G$ qui normalise $\Afr$. Il faudrait alors conna\^itre l'action de $\Pi_\Afr$
sur $\VV^{\lambda (\Afr )}$. Ce probl\`eme  devrait  \^etre r\'esolu en se donnant $\pi$ par un type
simple \'etendu. Jusqu'\`a pr\'esent, les types simples \'etendus ne sont construits que pour les
repr\'esentations supercuspidales (cf. [BK]{\S}6) et les repr\'esentations de la s\'erie discr\`ete de
niveau $0$ ([BH2], [SZ]).  
\end{remark2}

\section{R\'esolutions et pseudo-coefficients}

 En gardant les notations du {\S}3, on fixe un type simple $(J,\lambda
 )$ ainsi qu'une repr\'esentation $(\pi ,\VV )\in {\mathcal
   R}_{(J,\lambda )}$. On lui a associ\'e un $G$-complexe simplicial
 $X_\pi = X[L]$ muni d'un syst\`eme de coefficients $G$-\'equivariant
 ${\mathcal C}_\Theta (\pi )$. On peut alors former, pour $q= 0,
 ...,d-1$, $d={\rm dim}(X_\pi )$, les espaces $C_{q}^{\rm or}(X_\pi
 ,{\mathcal C}_\Theta (\pi ))$ de cha\^ines orient\'ees
 de $X_\pi$ \`a coefficients dans ${\mathcal C}_\Theta (\pi )$ ; ce
 sont des $G$-modules lisses. Ils sont munis d'applications bords
 $G$-\'equivariantes et  forment ainsi un complexe de cha\^ines
 naturellement augment\'e sur la repr\'esentation $\VV$ :
\begin{equation}
C_{d}^{\rm or}(X_\pi
 ,{\mathcal C}_\Theta (\pi ))
\stackrel{\partial}{\lra}  C_{d-1}^{\rm or}(X_\pi
 ,{\mathcal C}_\Theta (\pi ))
\stackrel{\partial}{\lra}
\cdots 
\stackrel{\partial}{\lra}
C_{0}^{\rm or}(X_\pi ,{\mathcal C}_\Theta (\pi ))
\stackrel{\epsilon}{\lra} \VV 
\lra 0
\end{equation}

\begin{proposition2} Le complexe augment\'e (1) est un complexe dans
  la cat\'egorie ${\mathcal R}_{(J,\lambda )}(G)$. 
\end{proposition2}

\noi {\bf Conjecture 1} {\it  Le complexe augment\'e (1) est une r\'esolution de
  $\VV$ dans la cat\'egorie ${\mathcal R}_{(J,\lambda )}(G)$.}
\medskip

Cette conjecture est une cons\'equence de la conjecture suivante :
\medskip

\noi {\bf Conjecture 2}. {\it Consid\'erons le complexe obtenu de (1) en
appliquant le foncteur (l'\'equivalence de cat\'egories) 
${\mathcal R}_{(J,\lambda )}(G) \lra {\mathcal
  H}_{\rm Spher}(G,\lambda_{\max})-{\rm Mod}$, ${\mathcal W}\mapsto
 {\mathcal W}^{\lambda_{\rm max}}$ :
 
\begin{equation}
C_{d}^{\rm or}(X_\pi
 ,{\mathcal C}_\Theta (\pi ))^{\lambda_{\rm max}}
\stackrel{\partial}{\lra}  C_{d-1}^{\rm or}(X_\pi
 ,{\mathcal C}_\Theta (\pi ))^{\lambda_{\rm max}}
\stackrel{\partial}{\lra}
\cdots 
\stackrel{\partial}{\lra}
C_{0}^{\rm or}(X_\pi ,{\mathcal C}_\Theta (\pi ))^{\lambda_{\rm max}}
\stackrel{\epsilon}{\lra} \VV^{\lambda_{\rm max}} 
\lra 0
\end{equation}

\noi Alors le complexe de cha\^ines augment\'e (2) est canoniquement
isomorphe au complexe de cha\^ines d'un appartement standard ${\mathcal
  A}_L$ de $X_L$ \`a coefficients constants dans $\VV^{\lambda_{\rm
    max}}$. Puisque ${\mathcal A}_L$ est contractile, ce dernier complexe est
exact.}
\medskip

Dans [BS], nous d\'emontrons la Conjecture 2, et donc la Conjecture 1 dans le cas
particulier suivant.

\begin{theorem2} {\it Avec les notations pr\'ec\'edentes, supposons que
  $(\pi ,\VV )$ est irr\'eductible et essentiellement de carr\'e
  int\'egrable. Alors la conjecture 2 est vraie.}
\end{theorem2}

Nous supposons dor\'enavant que la repr\'esentation $(\pi ,\VV )$ est
irr\'eductible et essentiellement de carr\'e int\'egrable. 
\medskip

 Soit $Z$ le centre de $G$. Fixons une mesure de Haar $\mu_{G/Z}$ sur
 $G/Z$. 
\smallskip

Pour chaque
simplexe $\sigma$ de $X_\pi$, on note $G_\sigma$ le stabilisateur
global de $\sigma$ dans $G$, et $\lambda_\sigma$ la repr\'esentation
irr\'eductible de $G_\sigma$ dans $\VV_\sigma $ induite par le syst\`eme de
coefficients ${\mathcal C}_\Theta (\pi )$ ; on note $\tau^\VV_\sigma$
le caract\`ere de $\lambda_\sigma$, que l'on \'etend par $0$ \`a $G$ en une
fonction localement constante \`a support compact modulo le centre ;
on note $\epsilon_\sigma~: \ G_\sigma \lra \{ \pm 1\}$ le caract\`ere
de $G_\sigma$ d\'efini comme suit : si $g\in G_\sigma$,
$\epsilon_\sigma (g)$ est la signature de la permutation des sommets
de $\sigma$ induite par $g$ ; enfin on \'etend la fonction
$\epsilon_\sigma$  \`a $G$ en une
fonction localement constante \`a support compact modulo le centre.

 Pour $q=0,...,d={\rm dim}\, X_\pi$, fixons un ensemble ${\mathcal
   F}_q$ de repr\'esentants des orbites de $G$ dans les $q$-simplexes
 de $X_\pi$. A la suite de  Kottwitz [Ko] et Schneider et Stuhler
 [SS], on attache \`a $(\pi ,\VV )$ une fonction dite {\it d'Euler-Poincar\'e}
  par la formule :
$$
f_{\rm EP}^\VV := \sum_{q=0}^d \sum_{\sigma \in {\mathcal F}_q} (-1)^q
\mu_{G/Z}(G_\sigma /Z )^{-1}\, {\bar
  \tau}_\sigma^\VV\ \epsilon_\sigma\ .
$$

En suivant une id\'ee originale de Kottwitz, reprise par Schneider
et Stuhler, on d\'emontre :

\begin{theorem2} La fonction $f_{\rm EP}^\VV$ est un
  pseudo-coefficient de $(\pi ,\VV )$.
\end{theorem2}

\section{La formule de caract\`ere}

 On fixe une repr\'esentation $(\pi ,\VV )$ de $G$ suppos\'ee
 irr\'eductibe et essentiellement de carr\'e int\'egrable. On note
 $\chi_\pi$ son caract\`ere d'Harish-Chandra. On garde
 les notations des sections pr\'ec\'edentes. Nous aurons besoin du
 r\'esultat suivant (d\^u  \`a Kazdhan [Ka] dans le cas d'un corps de base
 de caract\'eristique $0$ et d'un groupe r\'eductif \`a centre
 compact, et \`a Badulescu [Ba] dans le cas de notre groupe). 

\begin{theorem2} Soit $f_\pi$ un pseudo-coefficient et $\gamma\in G$ un
  \'el\'ement elliptique r\'egulier. Alors $\chi_\pi (\gamma )$ est donn\'e
  par l'int\'egrale orbitale convergente :
$$
\chi_\pi (\gamma ) = \int_{G/Z} f_{\pi} (g^{-1}\gamma^{-1}g)d\mu_{G/Z}
(\dot{g})\ .
$$
\end{theorem2}

 L'application de ce r\'esultat au pseudo-coefficient $f_{\rm EP}^\VV$
 donne une formule pour la valeur de
 $\chi_\pi$ en un \'el\'ement elliptique r\'egulier $\gamma$. Elle
 s'exprime de fa\c con \'el\'egante dans le cadre g\'eom\'etrique
 suivant. 

 Soit $\vert X_\pi \vert$ la r\'ealisation g\'eom\'etrique standard de
 $X_\pi$ ; c'est un espace topologique localement compact muni d'une
 action de $G$. L'ensemble de points fixes $\vert X_\pi \vert^\gamma$
 est compact. Il est muni de la  structure simpliciale naturelle  suivante : ses
 simplexes sont les intersections non vides  $\sigma (\gamma ):=\sigma
 \cap \vert X_\pi \vert^\gamma$, pour un simplexe $\sigma$ de $X_\pi$
 globalement fixe par $\gamma$. En fait $\sigma (\gamma )$ d\'etermine
 enti\`erement $\sigma$. 

 Avec ces notations, on a alors la formule de caract\`ere suivante. 

\begin{theorem2} La valeur du caract\`ere d'Harish-Chandra de $\pi$ en
  un \'el\'ement elliptique r\'egulier $\gamma$ est donn\'ee par
$$
\chi_\pi (\gamma )=\sum_{q=0}^{{\rm dim}\vert X_\pi\vert^\gamma}
  \sum_{\sigma (\gamma )\in \vert X_\pi\vert^\gamma_q} (-1)^q \  {\rm
    Tr}(\gamma ,\lambda_\sigma )\ ,
$$
\noi o\`u $\vert X_\pi\vert^\gamma_q$ d\'esigne l'ensemble des
$q$-simplexes de $\vert X_\pi\vert^\gamma$.
\end{theorem2}

 Cette derni\`ere formule ne peut en aucun cas être consid\'er\'ee
 comme {\it effective} en g\'en\'eral. En effet pour un \'el\'ement
 $\gamma$ elliptique r\'egulier quelconque, il n'existe aucune
 description connue de l'ensemble de points fixes  $\vert
 X_\pi\vert^\gamma$ (ni de l'ensemble  $\vert X \vert^\gamma$
 d'ailleurs). Il y a un cas cependant o\`u la formule se simplifie de
 fa\c con raisonnable, c'est celui o\`u $\gamma$ est minimal sur $F$
 au sens de [BK](1.4.14).

\begin{lemma2} Supposons $\gamma$ elliptique r\'egulier et  minimal
  sur $F$. Alors $\vert X\vert^\gamma = \vert X\vert^{K^\times}$, o\`u
  $K=F[\gamma ]$. En particulier $\vert X\vert^\gamma$ se r\'eduit \`a
  un point, image canonique de l'immeuble $X_K$ dans $X$. Dans un
  autre langage, $ \vert X\vert^\gamma =\{ x_\gamma \}$, o\`u
  $x_\gamma$ est l'isobarycentre du simplexe attach\'e \`a  l'unique
  ordre h\'er\'editaire $\Afr_\gamma$ de $A$ normalis\'e par
  $K^\times$ (ou de fa\c con \'equivalente par $\gamma$). 
\end{lemma2}

On en d\'eduit une {\it formule simple de caract\`ere} :

\begin{theorem2} Soit $\gamma$ un \'el\'ement elliptique r\'egulier de
  $G$, minimal sur $F$. Soit $\Afr_\gamma$ l'unique ordre
  h\'er\'editaire de $A$ normalis\'e par $F[\gamma ]^\times$, et soit
  $x_\gamma$ l'isobarycentre  du simplexe de $X$ correspondant \`a
  $\Afr_\gamma$. Si $x_\gamma \in \vert X_\pi \vert$, soit $\sigma_\gamma$
  l'unique simplexe de $X_\pi$ dont l'int\'erieur contient
  $x_\gamma$. Alors :
$$
\chi_\pi (\gamma )=\left\{
\begin{array}{ll}
{\rm Tr}(\gamma ,\lambda_{\sigma_\gamma}) & {\rm si}\ f(L/F)\vert
f(F[\gamma ]/F)\ {\rm et}\  e(L/F)\vert
e(F[\gamma ]/F),\\
0 & {\rm  sinon.}
\end{array}
\right.
$$
\end{theorem2}

Des cas particuliers de cette derni\`ere formule furent obtenus dans
[Br].

\bigskip

\centerline{\bf R\'ef\'erences}
\bigskip

[Ba] A.I. Badulescu, {\it Un r\'esultat de transfert et un
  r\'esultat d'int\'egrabilit\'e locale des caract\`eres en caract\'eristique
  non nulle},  J. Reine Angew. Math. 565 (2003), 101--124. 

[Br] P. Broussous, {\it Transfert du pseudo-coefficient de
  Kottwitz et formules de caractère pour la série discrète de GL(N)
  sur un corps local},  Canad. J. Math. 66 (2014), no. 2, 241–283. 

[BH] C.J. Bushnell et G. Henniart, {\it Local tame lifting
  for GL(N). I. Simple characters},  Inst. Hautes Études
  Sci. Publ. Math. No. 83 (1996), 105–233.

[BH2] C.J. Bushnell et G. Henniart, {\it Explicit
  functorial correspondences for level zero representations of $p$-adic
  linear groups}, J. Number Theory 131 (2011), no. 2, 309--331.

[BK]  C.J. Bushnell et P.C. Kutzko, {\it Smooth
  representations of reductive $p$-adic groups: structure theory via
  types}, Proc. London Math. Soc. (3) 77 (1998), 582--634.

[BL] P. Broussous et B. Lemaire, {\it Building of ${\rm
    GL}(m,D)$ and centralizers},  Transform. Groups 7 (2002), no. 1,
  15–50.

[BS] P. Broussous et P. Schneider, {\it Simple characters
  and coefficient systems on the building}, arXiv:1402.2501, 2014.  

[Ka] D. Kzhdan, {\it Cuspidal geometry of p-adic groups},
 J. Analyse Math. 47 (1986), 1–36.

[Ko] R.E. Kottwitz, {\it Tamagawa numbers},  Ann. of
  Math. (2) 127 (1988), no. 3, 629--646. 

[SS]  P. Schneider et  U. Stuhler, {\it Representation
  theory and sheaves  on the Bruhat-Tits building}, 
 Inst. Hautes Études Sci. Publ. Math. No. 85 (1997), 97--191. 

[SZ]  A.J. Silberger et E.-W. Zink,{\it  An explicit
  matching theorem for  level zero discrete series of unit groups of
  $p$-adic  simple algebras},  J. Reine Angew. Math. 585 (2005), 173--235.

\end{document}